\documentclass{article}[12pt]
\usepackage[a4paper, total={6in, 8in}]{geometry}

\usepackage{amsmath, amsthm, amssymb, amsfonts, mathtools, bbm}                                                                                                                                  
\usepackage{enumerate}
\usepackage{mathrsfs}

\newcommand{\norm}[1]{\lVert#1\rVert}

\newtheorem{theorem}{Theorem}
\newtheorem{lemma}{Lemma}
\newtheorem{proposition}{Proposition}

\theoremstyle{remark}

\title{A remark on the Beurling-Malliavin theorem in several variables}
\author{Alex Bergman}
\date{\today}

\begin{document}

\maketitle

\begin{abstract}
    We use a lifting trick to show that the Beurling-Malliavin multiplier theorem extends to radial functions in several variables in a straightforward way. This simplifies an argument of Vasilyev and also answers a question of Vasilyev on the Cartwright version of the theorem.
\end{abstract}

The first Beurling-Malliavin theorem on weights which can be minorized by functions with Fourier support contained in an arbitrarily small interval is one of the most celebrated results of 20th century one dimensional harmonic analysis. For a detailed discussion the reader may consult part two of the book \cite{MR1303780} or the survey paper \cite{MR2241422}. We recall one variant of the result below. We use the following normalization for the Fourier transform
\begin{equation*}
    \widehat{f}(\xi) = \int_{\mathbb{R}^{n}} f(x)e^{-ix\xi}dx.
\end{equation*}

\begin{theorem}[Beurling-Malliavin \cite{MR147848}]\label{thm:BM1}
    Let $\Omega : \mathbb{R} \to [0,\infty)$ be a Lipschitz continuous function satisfying
    \begin{equation*}
        \int_{-\infty}^{\infty} \frac{\Omega (x) }{1+x^{2}}dx < \infty.
    \end{equation*}
    Then for each $\epsilon > 0$ there exists a non-trivial function $f \in L^{2}(\mathbb{R})$ with $\text{supp}(\widehat{f}) \subset (-\epsilon, \epsilon)$, such that
    \begin{equation*}
        \lvert f(x) \rvert \leq \exp(-\Omega(x)), \; \text{ for } x \in \mathbb{R}.
    \end{equation*}
\end{theorem}
The first Beurling-Malliavin theorem has many striking applications. Let us mention the second Beurling-Malliavin theorem on the radius of completeness of systems of exponentials \cite{MR209758}, the work of Eremenko and Novikov on oscillations of high-pass signals \cite{MR2038065, MR2048457}, and the fractal uncertainty principle of Bourgain and Dyatlov \cite{MR3779959}. Partly due to its many applications and partly because of the beauty of the result there has been some recent effort in understanding multi-dimensional generalizations of this theorem, see \cite{MR4927737, MR4884567}. In this note we show that the results in \cite{MR4884567} follow from the one-dimensional theorem and a lifting trick and the Paley-Wiener theorem adapted to several variables. In addition, we show that the same method also applies to the Cartwright version of the theorem.

\section{The radial Beurling-Malliavin theorem}

    The Paley-Wiener theorem in several variables asserts the equality of the following two classes
    \begin{enumerate}[(i)]
        \item Functions $f \in L^{2}(\mathbb{R}^{n})$, such that $\text{supp}(\widehat{f}) \subset B(0,\sigma)$.
        \item Entire functions $f : \mathbb{C}^{n} \to \mathbb{C}$ satisfying
        \begin{equation*}
            \lvert f(x+iy) \rvert \leq Ae^{\sigma \norm{y}}, \; \text{ for all } x,y \in \mathbb{R}^{n} \; \text{(} \norm{\cdot} \text{ is the usual Euclidean norm}), 
        \end{equation*}
        for some constant $A$ (possibly depending on $f$) and whose restriction to $\mathbb{R}^{n}$ belongs to $L^{2}(\mathbb{R}^{n})$.
    \end{enumerate}
    See, for example, Hörmander \cite{MR1065993}, Chapter 7.3. We shall need the following simple lifting lemma.
    \begin{lemma}\label{lemma:lemma1}
        Let $f_{0} : \mathbb{C} \to \mathbb{C}$ be an even entire function with power series expansion
        \begin{equation*}
            f_{0}(z) = \sum_{m=0}^{\infty}a_{m}\zeta^{2m}, \; \text{ for } \zeta \in \mathbb{C},
        \end{equation*}
        satisfying
        \begin{equation*}
            \lvert f(a+ib) \rvert \leq A e^{\sigma \lvert b \rvert}, \; \text{ for all } a,b \in \mathbb{R}.
        \end{equation*}
        Then the entire function $f : \mathbb{C}^{n} \to \mathbb{C}$ given by
        \begin{equation*}
            f(z) = \sum_{m=0}^{\infty}a_{m}(z_{1}^{2}+...+z_{n}^{2})^{m}, \; \text{ for } z \in \mathbb{C}^{n},
        \end{equation*}
        satisfies
        \begin{equation*}
            \lvert f(x+iy) \rvert \leq Ae^{\sigma \norm{y}}, \; \text{ for } x,y \in \mathbb{R}^{n}.
        \end{equation*}
        \begin{proof}
            Fix $z = (z_{1},...,z_{n}) \in \mathbb{C}^{n}$, $z = x+iy$. Let $\zeta = a+ib \in \mathbb{C}$ be a complex number such that
            \begin{equation*}
                \zeta^{2} = z_{1}^{2}+...+z_{n}^{2}.
            \end{equation*}
            Then
            \begin{equation*}
                \lvert f(x+iy)\rvert = \lvert f_{0}(\zeta)\rvert \leq Ae^{\sigma \lvert b \rvert}.    
            \end{equation*}
            If $\lvert b \rvert \leq \norm{y}$ we are done. Suppose, for a contradiction, that $\norm{y} < \lvert b \rvert$. We have
            \begin{equation*}
                \begin{cases}
                    a^{2}-b^{2} = \norm{x}^{2} - \norm{y}^{2}, \\
                    ab = \langle x, y\rangle.
                \end{cases}
            \end{equation*}
            Thus
            \begin{equation*}
                a^{2}-\norm{x}^{2} = b^{2}-\norm{y}^{2} > 0.
            \end{equation*}
            Hence $\lvert a \rvert > \norm{x}$. It follows from this and the equation
            \begin{equation*}
                \lvert ab \rvert \leq \norm{x}\norm{y},
            \end{equation*}
            that
            \begin{equation*}
                \lvert b \rvert \leq \norm{y}.
            \end{equation*}
            this is a contradiction.
        \end{proof}
    \end{lemma}

    Let us begin with the radial version of the first Beurling-Malliavin theorem from \cite{MR4884567}.
    \begin{theorem}[Vasilyev \cite{MR4884567}]\label{thm:radial_BM1}
        Let $\Omega : \mathbb{R}^{n} \to [0,\infty)$ be a radial Lipschitz continuous function satisfying
    \begin{equation*}
        \int_{\mathbb{R}^{n}} \frac{\Omega (x)dx}{(1+\norm{x}^{2})^{\frac{n+1}{2}}} < \infty.
    \end{equation*}
    Then for each $\epsilon > 0$ there exists a non-trivial function $f \in L^{2}(\mathbb{R})$ with $\text{supp}(\widehat{f}) \subset B(0,\epsilon)$, such that
    \begin{equation*}
        \lvert f(x) \rvert \leq e^{-\Omega(x)}, \; \text{ for } x \in \mathbb{R}^{n}.
    \end{equation*}
    \begin{proof}
        Fix $\epsilon > 0$. Let $\Omega_{0} : \mathbb{R} \to [0,\infty)$ denote the radial part of $\Omega$ extended to be an even function. Then $\Omega_{0}$ is Lipschitz continuous and
        \begin{equation*}
            \infty >\int_{\mathbb{R}^{n}} \frac{\Omega(x)dx}{(1+\norm{x}^{2})^{\frac{n+1}{2}}} = \frac{\lvert S^{n-1}\rvert}{2}\int_{-\infty}^{\infty} \frac{\Omega_{0}(r)\lvert r\rvert^{n-1}dr}{(1+r^{2})^{\frac{n+1}{2}}}.
        \end{equation*}
        Thus by the one dimensional Beurling-Malliavin theorem there exists a function $f_{0} \in L^{2}(\mathbb{R})$ with $\text{supp}(\widehat{f_{0}}) \subset (-\epsilon, \epsilon)$ and
        \begin{equation*}
            \lvert f_{0}(x) \rvert \leq e^{-\Omega_{0}(x)}, \; \text{ for all } x \in \mathbb{R}.
        \end{equation*}
        Since $\Omega_{0}$ is even $f_{0}$ may be chosen to be even (otherwise replace $f_{0}$ by $2^{-1}(f_{0}(x)+f_{0}(-x))$). Moreover, we may also suppose that
        \begin{equation*}
            \int_{-\infty}^{\infty} \lvert f_{0}(x) \rvert^{2} \lvert x \rvert^{n-1}dx < \infty.    
        \end{equation*}
        Indeed, we may choose a function in $(-\epsilon/2,\epsilon/2)$ minorizing $\exp(-\Omega)$ and multiply it by a suitable Schwartz function with small spectrum to obtain the required integrability. Since $f_{0}$ is even it has the power series expansion
        \begin{equation*}
            f_{0}(z) = \sum_{m=0}^{\infty} a_{m}z^{2m}.
        \end{equation*}
        Define the $\mathbb{C}^{n}$ entire function $f$ by
        \begin{equation*}
            f(z) = \sum_{m=0}^{\infty} a_{m}(z_{1}^{2}+...+z_{n}^{2})^{m}, \; \text{ for } z\in \mathbb{C}^{n}.
        \end{equation*}
        Then
        \begin{equation*}
            \lvert f(x) \rvert = \lvert f_{0}(\norm{x}) \rvert \leq e^{-\Omega_{0}(\norm{x})} = e^{-\Omega(x)}, \; \text{ for all } x \in \mathbb{R}^{n}.
        \end{equation*}
        Moreover,
        \begin{equation*}
            \int_{\mathbb{R}^{n}} \lvert f(x)\rvert^{2}dx = \lvert S^{n-1} \rvert \int_{0}^{\infty} \lvert f_{0}(r)\rvert^{2} r^{n-1}dr < \infty.
        \end{equation*}
        It remains to prove that $\widehat{f}$ is supported in $B(0,\epsilon)$. For this we use the Paley-Wiener theorem in several variables and Lemma \ref{lemma:lemma1}. Indeed, since $f_{0}$ satisfies
        \begin{equation*}
            \lvert f_{0}(a+ib) \rvert \leq C(f_{0}) e^{\epsilon \lvert b \rvert}, \text{ for all } a,b \in \mathbb{R},
        \end{equation*}
        It follows from Lemma \ref{lemma:lemma1} that
        \begin{equation*}
            \lvert f(x+iy)\rvert \leq C(f_{0}) e^{\epsilon \norm{y}}, \; \text{ for all } x,y \in \mathbb{R}^{n}.
        \end{equation*}
        Hence by the Paley-Wiener theorem $\text{supp}(\widehat{f}) \subset B(0,\epsilon)$.
    \end{proof}
    \end{theorem}

    \subsection{The multiplier theorem}

    Now let us consider the other variant of the first Beurling-Malliavin theorem concerning weights of the form $\omega(x) = \lvert f(x) \rvert^{-1}$ where $f$ belongs to the Cartwright class. To motivate the definition of the Cartwright class in several variables we show that the integrability condition on $\Omega$ cannot be dropped.

    \begin{proposition}
        Let $\Omega : \mathbb{R} \to [0,\infty)$ be a radial measurable function and suppose that there exists a non-trivial $f \in L^{2}(\mathbb{R}^{n})$ with $\text{supp}(\widehat{f}) \subset B(0,\sigma)$, such that
        \begin{equation*}
            \lvert f(x) \rvert \leq e^{-\Omega(x)}, \; \text{ for almost every } x \in \mathbb{R}^{n}.
        \end{equation*}
        Then
        \begin{equation*}
            \int_{\mathbb{R}^{n}} \frac{\Omega(x)dx}{(1+\norm{x}^{2})^{\frac{n+1}{2}}}<\infty.
        \end{equation*}
        \begin{proof}
            Replacing $f$ by $h = \lvert f \rvert^{2}$ and $\Omega$ by $2\Omega$ we obtain
            \begin{equation*}
                h(x) \leq e^{-2\Omega(x)}, \; \text{ for almost every } x \in \mathbb{R}^{n}.
            \end{equation*}
            with $h \in L^{1}(\mathbb{R}^{n})$ nonnegative and $\text{supp}(\widehat{h}) \subset B(0,2\sigma)$. Since $\Omega$ is radial we claim that $h$ can be chosen to be radial. Indeed, let $O(n)$ denote the group of orthogonal $n\times x$ matrices. Let $dU$ be the right-invariant Haar measure of the group normalized to be a probability measure. Since $O(n)$ is compact the measure is actually bi-invariant. By assumption
            \begin{equation*}
                \lvert h(Ux) \rvert \leq e^{-2\Omega(x)}, \; \text{ for almost every } x \in \mathbb{R}^{n} \text{ and } U \in O(n).
            \end{equation*}
            Integrating with respect to the Haar measure gives
            \begin{equation*}
                \int_{O(n)} \lvert h(Ux) \rvert dU \leq e^{-2\Omega(x)}, \; \text{ for almost every } x \in \mathbb{R}^{n}
            \end{equation*}
            Thus the symmetrized function
            \begin{equation*}
                \Tilde{h}(z) = \int_{O(n)} h(Uz) dU,
            \end{equation*}
            minorizes $\exp(-2\Omega)$. It is non-trivial since $h$ is non-trivial and $h \geq 0$. Since it is summable and of exponential type at most the exponential type of $h$ the claim follows. Let $h_{0}$ denote the radial component of $h$ extended to be an even function. Then $h_{0}$ is an entire function of exponential type and by considering polar coordinates we see that $x^{n-1}h_{0}$ is summable on the real-line. Hence, $h_{0}$ is a Cartwright class function. Thus
            \begin{equation*}
                -\int_{\mathbb{R}^{n}} \frac{\Omega(x)dx}{(1+\norm{x}^{2})^{\frac{n+1}{2}}} \geq \frac{1}{2} \int_{\mathbb{R}^{n}} \frac{\log \lvert h(x) \rvert}{(1+\norm{x}^{2})^{\frac{n+1}{2}}} = \frac{\lvert S^{n-1} \rvert}{2} \int_{0}^{\infty} \frac{\log \lvert h_{0} (r)\rvert}{(1+r^{2})^{\frac{n+1}{2}}}r^{n-1} > -\infty.
        \end{equation*}
        \end{proof}
    \end{proposition}

    Thus it is natural to introduce the Cartwright class of exponential type $\sigma$, $C_{\sigma}(\mathbb{C}^{n})$ in several complex variables as the set of entire functions $f : \mathbb{C}^{n} \to \mathbb{C}$ of satisfying
    \begin{equation*}
       \lvert f(x+iy) \rvert \leq Ae^{\sigma \norm{y}}, \; \text{ for all } x,y \in \mathbb{R}^{n}, 
    \end{equation*}
    and
    \begin{equation*}
            \int_{\mathbb{R}^{n}} \frac{\lvert \log \lvert f(x) \rvert \rvert dx}{(1+\norm{x}^{2})^{\frac{n+1}{2}}} < \infty.
    \end{equation*}
    The Cartwright class is
    \begin{equation*}
        C(\mathbb{C}^{n}) = \bigcup_{a > 0} C_{a}(\mathbb{C}^{n}).
    \end{equation*}
    As in one variable the (radial) Paley-Wiener spaces are contained in the Cartwright classes.
    \begin{proposition}
        Let $f \in L^{2}(\mathbb{R}^{n})$ with $\text{supp}(\widehat{f}) \subset B(0,\sigma)$ be radial and not identically zero. Then $f \in C_{\sigma}(\mathbb{C}^{n})$.
        \begin{proof}
            By the Paley-Wiener theorem everything is clear except the logarithmic integrability of $f$. Since $f \in L^{2}(\mathbb{R}^{n})$ it suffices to prove that
            \begin{equation*}
                \int_{\mathbb{R}^{n}} \frac{\log \lvert f(x)\rvert dx}{(1+\norm{x}^{2})^{\frac{n+1}{2}}} > -\infty.
            \end{equation*}
           The radial part, $f_{0}$, extended to be an even function is again an entire function of exponential type $\sigma$ which belongs to the one-dimensional Cartwright class. Hence,
            \begin{equation*}
                \int_{\mathbb{R}^{n}} \frac{\log \lvert f(x)\rvert dx}{(1+\norm{x}^{2})^{\frac{n+1}{2}}} = \lvert S^{n-1} \rvert \int_{0}^{\infty} \frac{\log \lvert f_{0}(r)\rvert r^{n-1}dr}{(1+r^{2})^{\frac{n+1}{2}}}> -\infty.
            \end{equation*}
        \end{proof}
    \end{proposition}
    
    We have the following radial variant in higher dimensions which was asked for in \cite{MR4884567}.

    \begin{theorem}
        Let $f \in C(\mathbb{C}^{n})$ be radial. Then for every $\epsilon > 0$ there exists $G \in C_{\epsilon}(\mathbb{C}^{n})$, such that
        \begin{equation*}
            \lvert G(x)f(x) \rvert \leq 1, \text{ for all } x \in \mathbb{R}^{n}.
        \end{equation*}
        \begin{proof}
            Since $f$ is an entire function invariant under the orthogonal group it has the expansion
            \begin{equation*}
                f(z) = \sum_{m=0}^{\infty}a_{m}(z_{1}^{2}+...+z_{n}^{2})^{m}, \; \text{ for } z \in \mathbb{C}^{n},
            \end{equation*}
            see Lemma 2.11. in Chapter 4 of \cite{MR304972}. Define the even entire function of exponential type $f_{0}$ by
            \begin{equation*}
                f_{0}(\zeta) = \sum_{m=0}^{\infty}a_{m}\zeta^{2m}, \; \text{ for } \zeta \in \mathbb{C}.
            \end{equation*}
            Applying polar coordinates we see that $f_{0}$ belongs to the usual one dimensional Cartwright class of the same exponential type as $f$. Thus by the one dimensional Beurling-Malliavin theorem (for weights of the form $\omega(x) = \lvert h(x) \rvert^{-1}$ with $h$ in the Cartwright class) there exists, $G_{0}$ in the one dimensional Cartwright class of exponential type at most $\epsilon$, such that
            \begin{equation*}
                \lvert G_{0}(x) f_{0}(x)\rvert \leq 1, \; \text{ for all } x\in \mathbb{R}.
            \end{equation*}
            Moreover, as before, without loss of generality we may assume $G_{0}$ is even. Hence, it has power series expansion
            \begin{equation*}
                G_{0}(\zeta) = \sum_{m=0}^{\infty} b_{m}\zeta^{2m}, \text{ for } \zeta \in \mathbb{C}.
            \end{equation*}
            Define
            \begin{equation*}
                G(z) = \sum_{m=0}^{\infty} b_{m}(z_{1}^{2}+...+z_{n}^{2})^{m}, \text{ for } z \in \mathbb{C}^{n}.
            \end{equation*}
            As before, $G$ is an entire function,
            \begin{equation*}
                \lvert G(z) \rvert \leq A e^{\epsilon \norm{y}},
            \end{equation*}
            and
            \begin{equation*}
                \lvert G(x)f(x) \rvert = \lvert G_{0}(\norm{x})f_{0}(\norm{x}) \rvert \leq 1, \; \text{ for all } x \in \mathbb{R}^{n}.
            \end{equation*}
            The required logarithmic integrability of $G$ follows from the corresponding integrability condition for $G_{0}$.
        \end{proof}
    \end{theorem}

\section{The non-radial case}

    Let us briefly discuss what can be said for the non-radial case. Let $\omega = \exp(-\Omega)$ be a weight function, $\Omega : \mathbb{R}^{n} \to [0,\infty)$. Our goal will be to construct a radial Lipschitz function $\Omega_{1}$ majorizing $\Omega$ and satisyfing
    \begin{equation*}
        \int_{\mathbb{R}^{n}} \frac{\Omega_{1}(x)dx}{(1+\norm{x}^{2})^{\frac{n+1}{2}}} < \infty.
    \end{equation*}
    Then for every $\epsilon > 0$, by what was proved for the radial case, there exists a a nontrivial function $f \in L^{2}(\mathbb{R}^{n})$ with Fourier support contained in $B(0,\epsilon)$, such that
    \begin{equation*}
        \lvert f(x) \rvert \leq e^{-\Omega_{1}(x)} \leq e^{-\Omega(x)}, \; \text{ for all } x \in \mathbb{R}^{n}.
    \end{equation*}
    The results of \cite{MR4927737} are not attainable (in an obvious way) through this method since there the weight is first averaged over lines
    \begin{equation*}
        G(x) = \int_{1/2}^{2}\Omega(sx)ds, \; \text{ for } x \in \mathbb{R}^{n},
    \end{equation*}
    and then made radial
    \begin{equation*}
        G^{*}(r) = \sup_{\norm{x}=r} G(x), \; \text{ for } r \geq 0.
    \end{equation*}
    Now the relationship between $G^{*}$ and $\Omega$ is somewhat nontrivial. Where as, in \cite{MR4884567} $\Omega_{1}$ is constructed by taking the supremum of $\Omega$ on dyadic annuli $B(0,2^{j+1})\setminus B(0,2^{j})$. Thus in the second case there is no averaging and it is clear that $\Omega_{1}$ dominates $\Omega$.


\bibliographystyle{abbrv}
\bibliography{citations}

\begin{thebibliography}{10}

\bibitem{MR147848}
A.~Beurling and P.~Malliavin.
\newblock On {F}ourier transforms of measures with compact support.
\newblock {\em Acta Math.}, 107:291--309, 1962.

\bibitem{MR209758}
A.~Beurling and P.~Malliavin.
\newblock On the closure of characters and the zeros of entire functions.
\newblock {\em Acta Math.}, 118:79--93, 1967.

\bibitem{MR3779959}
J.~Bourgain and S.~Dyatlov.
\newblock Spectral gaps without the pressure condition.
\newblock {\em Ann. of Math. (2)}, 187(3):825--867, 2018.

\bibitem{MR4927737}
A.~Cohen.
\newblock Fractal uncertainty in higher dimensions.
\newblock {\em Ann. of Math. (2)}, 202(1):265--307, 2025.

\bibitem{MR2038065}
A.~Eremenko and D.~Novikov.
\newblock Oscillation of {F}ourier integrals with a spectral gap.
\newblock {\em J. Math. Pures Appl. (9)}, 83(3):313--365, 2004.

\bibitem{MR2048457}
A.~Eremenko and D.~Novikov.
\newblock Oscillation of functions with a spectral gap.
\newblock {\em Proc. Natl. Acad. Sci. USA}, 101(16):5872--5873, 2004.

\bibitem{MR1303780}
V.~Havin and B.~J\"oricke.
\newblock {\em The uncertainty principle in harmonic analysis}, volume~28 of
  {\em Ergebnisse der Mathematik und ihrer Grenzgebiete (3) [Results in
  Mathematics and Related Areas (3)]}.
\newblock Springer-Verlag, Berlin, 1994.

\bibitem{MR1065993}
L.~H\"ormander.
\newblock {\em The analysis of linear partial differential operators. {I}},
  volume 256 of {\em Grundlehren der mathematischen Wissenschaften [Fundamental
  Principles of Mathematical Sciences]}.
\newblock Springer-Verlag, Berlin, second edition, 1990.
\newblock Distribution theory and Fourier analysis.

\bibitem{MR2241422}
D.~Mashregi, F.~L. Nazarov, and V.~P. Khavin.
\newblock The {B}eurling-{M}alliavin multiplier theorem: the seventh proof.
\newblock {\em Algebra i Analiz}, 17(5):3--68, 2005.

\bibitem{MR304972}
E.~M. Stein and G.~Weiss.
\newblock {\em Introduction to {F}ourier analysis on {E}uclidean spaces},
  volume No. 32 of {\em Princeton Mathematical Series}.
\newblock Princeton University Press, Princeton, NJ, 1971.

\bibitem{MR4884567}
I.~Vasilyev.
\newblock The {B}eurling and {M}alliavin theorem in several dimensions.
\newblock {\em Math. Ann.}, 391(4):6057--6072, 2025.

\end{thebibliography}

\end{document}